\documentclass{compositio}
\usepackage{mathrsfs}
\usepackage{amsmath}
\usepackage{amssymb}
\usepackage{verbatim}
\usepackage{xcolor}
\usepackage{graphicx}
\usepackage[all]{xy}
\numberwithin{equation}{section}

\theoremstyle{plain}
\newtheorem{theorem}{Theorem}[section]
\newtheorem{proposition}{Proposition}[section]

\theoremstyle{definition}
\newtheorem{definition}{Definition}[section]

\theoremstyle{remark}
\newtheorem{rem}{Remark}[section]

\begin{document}

\title{The Kawamata--Viehweg--Nadel-type vanishing theorem and the asymptotic multiplier ideal sheaf}

\author{Jingcao Wu}
\email{wujincao@shufe.edu.cn}
\address{School of Mathematics, Shanghai University of Finance and Economics, Shanghai 200433, People's Republic of China}

\thanks{This research was supported by China Postdoctoral Science Foundation, grant 2019M661328.}

\begin{abstract}
In this paper we establish a Nadel-type vanishing theorem and a Kawamata--Viehweg-type vanishing theorem concerning the asymptotic multiplier ideal sheaf on a compact K\"{a}hler manifold $X$. After that, we provide a relative variant.
\end{abstract}

\maketitle

\section{Introduction}
\label{sec:introduction}

The celebrated Nadel vanishing theorem says that
\begin{theorem}[(c.f. \cite{Nad90,Dem93})]\label{t11}
Let $X$ be a projective manifold of dimension $n$, and let $L$ be a holomorphic line bundle on $X$. Fix a K\"{a}hler metric $\omega$ on $X$. Assume that $L$ is equipped with a singular Hermitian metric $\varphi$ with $i\Theta_{L,\varphi}\geqslant\varepsilon\omega$ for some $\varepsilon>0$. Then
\[
H^{q}(X,K_{X}\otimes L\otimes\mathscr{I}(\varphi))=0
\]
for $q>0$.
\end{theorem}

Here $\mathscr{I}(\varphi)$ refers to the multiplier ideal sheaf \cite{Nad90} associated to $\varphi$. Note that in Theorem \ref{t11}, $L$ is by definition a big line bundle \cite{Dem12}. So we will also directly say that $(L,\varphi)$ is a big line bundle. 

This theorem can be seen as the analytic counterpart of the Kawamata--Viehweg vanishing theorem \cite{Kaw82,Vie82} in algebraic geometry. The later one is stated as follows:
\begin{theorem}[(c.f. \cite{Kaw82,Vie82})]\label{t12}
Let $X$ be a projective manifold of dimension $n$, and let $L$ be a holomorphic line bundle on $X$. Assume that $L$ is numerically equivalent to $B+\Delta$, where $B$ is a nef and big $\mathbb{Q}$-divisor, and $\Delta=\sum a_{i}\Delta_{i}$ is a $\mathbb{Q}$-divisor with simple normal crossing support and fractional coefficients:
\[
0\leqslant a_{i}<1 \textrm{ for all i}.
\]
Then we have
\[
   H^{q}(X,K_{X}\otimes L)=0
\]
for $q>0$.
\end{theorem} 

Both of them are fundamental tools in complex geometry. However, one of the limitations of Theorem \ref{t11} (resp. Theorem \ref{t12}) is that $L$ (resp. $B$) should be big. It is then asked to generalise these two vanishing theorems. One could refer to \cite{Cao14,Eno93,Mat14,Mat15a,Mat15b,Mat18,Wu17} and the references therein for several generalisations. In practice, we find that the lower bound of the order $q$ such that $H^{q}(X,K_{X}\otimes L\otimes\mathscr{I}(\varphi))=0$ usually depends on the numerical dimension $\textrm{nd}(L)$ or Iitaka dimension $\kappa(L)$ \cite{Laz04} of $L$. In particular, in Theorem \ref{t11}, $\kappa(L)=n$, and the lower bound is just $0=n-\kappa(L)$.

In this paper, we present the following Nadel-type vanishing theorem concerning the asymptotic multiplier ideal sheaf $\mathscr{I}(\|L\|)$ (see Sect.\ref{sec:asymptotic}).
\begin{theorem}\label{t13}
Let $X$ be a compact K\"{a}hler manifold of dimension $n$, and let $L$ be a pseudo-effective line bundle on $X$. Then we have
\[
   H^{q}(X,K_{X}\otimes L\otimes\mathscr{I}(\|L\|))=0
\]
for $q>n-\kappa(L)$.
\end{theorem}

As for the algebraic counterpart, we have the following generalisation of the Kawamata--Viehweg vanishing theorem. 
\begin{theorem}\label{t14}
Let $X$ be a compact K\"{a}hler manifold of dimension $n$, and let $L$ be a holomorphic line bundle on $X$. Assume that $L$ is numerically equivalent to $B+\Delta$, where $B$ is a nef and abundant $\mathbb{Q}$-divisor, and $\Delta=\sum a_{i}\Delta_{i}$ is a $\mathbb{Q}$-divisor with simple normal crossing support and fractional coefficients:
\[
0\leqslant a_{i}<1 \textrm{ for all i}.
\]
Then we have
\[
   H^{q}(X,K_{X}\otimes L)=0
\]
for $q>n-\textrm{nd}(B)$.
\end{theorem}
We will make a brief introduction for the abundant divisor in Sect. \ref{sec:abundant}. In particular, the asymptotic multiplier ideal sheaf precisely characterises this notion. (See Theorem \ref{t21}.)

\begin{rem}\label{r11}
In the original manuscript, Theorem \ref{t13} was established on a projective manifold. The referees kindly pointed out the possibilities to extend the whole things to the K\"{a}hler case, which leads to the current form. We sincerely appreciate it.
\end{rem}

Eventually, we present a relative version of Theorem \ref{t13}. 
\begin{theorem}\label{t15}
Let $f:X\rightarrow Y$ be a surjective morphism of projective varieties, with $X$ non-singular. Let $L$ be a pseudo-effective line bundle on $X$. Let $l$ be the dimension of a general fibre $F$ of $f$. Then
\[
R^{q}f_{\ast}(K_{X}\otimes L\otimes\mathscr{I}(f,\|L\|))=0
\]
for $q>l-\kappa(L,f)$. Here $\mathscr{I}(f,\|L\|)$ is the relative version of the asymptotic multiplier ideal sheaf (see Sect.\ref{sec:asymptotic}) and $\kappa(L,f)$ is the relative Iitaka dimension (see Sect.\ref{sec:iitaka}).
\end{theorem}
\begin{rem}\label{r12}
Note that Theorem \ref{t15} is not easily obtained by applying Theorem \ref{t13} on the general fibre. Moreover, when proving Theorem \ref{t13}, we use the Monge--Amp\`{e}re technique in \cite{DeP03}. It is difficult to extend it to the relative setting. Instead, we directly prove Theorem \ref{t15} by induction on the dimension. However, this method clearly fails for a non-projective setting. We hope this limitation will be overcome in the future.
\end{rem}

\begin{rem}\label{r13}
The relative version of Theorem \ref{t14} is similar, and we will omit it here.
\end{rem}

This paper is organised as follows. We first recall some background materials, including the asymptotic multiplier ideal sheaf, abundant line bundle and so on. Then, we proceed to prove Theorems \ref{t13} and \ref{t14} in Sect.\ref{sec:absolute}. In the end, we prove Theorem \ref{t15}.

\begin{acknowledgements}
The author wants to thank Prof. Jixiang Fu for his suggestion and encouragement. Thanks are also due to the referees for pointing out the possibilities to extend the main theorem to the K\"{a}hler case, which brings new value to this paper.
\end{acknowledgements}

\section{Preliminary}
\label{sec:preliminary}
In this section we will introduce some basic materials.  

\subsection{The asymptotic multiplier ideal sheaf}
\label{sec:asymptotic}
This part is mostly collected from \cite{Laz04b}.

Let $X$ be a compact K\"{a}hler manifold. Recall that for an arbitrary ideal sheaf $\mathfrak{a}\subset\mathcal{O}_{X}$, the associated multiplier ideal sheaf is defined as follows: let $\mu:\tilde{X}\rightarrow X$ be a smooth modification such that $\mu^{\ast}\mathfrak{a}=\mathcal{O}_{\tilde{X}}(-E)$, where $E$ has simple normal crossing support. Then given a positive real number $c>0$ the multiplier ideal sheaf is defined as
\[
\mathscr{I}(c\cdot\mathfrak{a}):=\mu_{\ast}\mathcal{O}_{\tilde{X}}(K_{\tilde{X}/X}-\lfloor cE\rfloor).
\]
Here $K_{\tilde{X}/X}$ is the relative canonical bundle and $\lfloor E\rfloor$ means the round-down.

Let $L$ be a holomorphic line bundle on $X$ with $\kappa(L)\geqslant0$. Fix a positive real number $c>0$. For $k>0$ consider the complete linear series $|L^{k}|$, and form the multiplier ideal sheaf
\[
\mathscr{I}(\frac{c}{k}|L^{k}|)\subseteq\mathcal{O}_{X},
\]
where $\mathscr{I}(\frac{c}{k}|L^{k}|):=\mathscr{I}(\frac{c}{k}\cdot\mathfrak{a}_{k})$ with $\mathfrak{a}_{k}$ being the base-ideal of $|L^{k}|$. It is not hard to verify that for every integer $p\geqslant1$ one has the inclusion
\[
\mathscr{I}(\frac{c}{k}|L^{k}|)\subseteq\mathscr{I}(\frac{c}{pk}|L^{pk}|).
\]
Therefore the family of ideals
\[
\{\mathscr{I}(\frac{c}{k}|L^{k}|)\}_{(k\geqslant0)}
\]
has a unique maximal element from the ascending chain condition on ideals.

\begin{definition}\label{d21}
The asymptotic multiplier ideal sheaf associated to $c$ and $|L|$,
\[
\mathscr{I}(c\|L\|)
\]
is defined to be the unique maximal member among the family of ideals $\{\mathscr{I}(\frac{c}{k}|L^{k}|)\}$.
\end{definition}
Observe that 
\begin{equation}\label{c21}
\mathscr{I}(c\|L^{p}\|)=\mathscr{I}(cp\|L\|)
\end{equation}
for any integer $p\geqslant1$. (See Theorem 11.1.8 in \cite{Laz04b}, for example.) Then we could extend Definition \ref{d21} to a $\mathbb{Q}$-divisor $D$ as follows: take positive integer $p$ such that $\mathcal{O}_{X}(pD)$ is a line bundle. Then for any real number $c>0$,
\[
\mathscr{I}(c\|D\|):=\mathscr{I}(\frac{c}{p}\|pD\|).
\]
By (\ref{c21}) we see that $\mathscr{I}(c\|D\|)$ is independent of the choice of $p$, hence is well-defined.

Now we return to the line bundle, or $\mathbb{Z}$-divisor case. By definition, $\mathscr{I}(c\|L\|)=\mathscr{I}(\frac{c}{k}|L^{k}|)$ for some $k$. In this case, we say that $k$ computes $\mathscr{I}(c\|L\|)$. Let $\{u_{1},...,u_{m}\}$ be a basis of $H^{0}(X,L^{k})$, then the base-ideal of $|L^{k}|$ is just $\mathcal{I}(u_{1},...,u_{m})$. Let $\varphi=\log(|u_{1}|^{2}+\cdots+|u_{m}|^{2})$, which is a singular metric on $L^{k}$. We verify that
\[
\mathscr{I}(\frac{c}{k}|L^{k}|)=\mathscr{I}(\frac{c}{k}\varphi).
\]
Indeed, let $\mu:\tilde{X}\rightarrow X$ be the smooth modification such that $\mu^{\ast}\mathcal{I}(u_{1},...,u_{m})=\mathcal{O}_{\tilde{X}}(-E)$, where $E$ has simple normal crossing support. Then it is computed in \cite{Dem12} that
\[
\mathscr{I}(\frac{c}{k}\varphi)=\mu_{\ast}\mathcal{O}_{\tilde{X}}(K_{\tilde{X}/X}-\lfloor \frac{c}{k}E\rfloor).
\]
Obviously, this representation of $\mathscr{I}(\frac{c}{k}\varphi)$ coincides with the definition of $\mathscr{I}(\frac{c}{k}|L^{k}|)$. In summary, we have
\[
\mathscr{I}(c\|L\|)=\mathscr{I}(\frac{c}{k}\varphi),
\]
and $\frac{1}{k}\varphi$ is called the singular metric on $L$ associated to $\mathscr{I}(c\|L\|)$. Certainly the associated metric depends on the choice of $k$, hence is not unique.

Next, we introduce the relative variant. Let $f:X\rightarrow Y$ be a surjective morphism between compact K\"{a}hler manifolds, and let $L$ be a holomorphic line bundle on $X$ whose restriction to a general fibre of $f$ has non-negative Iitaka dimension. Then there is a naturally defined homomorphism
\[
\rho:f^{\ast}f_{\ast}(L)\rightarrow L.
\]
Let $\mu:\tilde{X}\rightarrow X$ be a smooth modification of $|L|$ with respect to $f$, having the property that the image of the induced homomorphism
\[
\mu^{\ast}\rho:\mu^{\ast}f^{\ast}f_{\ast}(L)\rightarrow \mu^{\ast}L
\]
is the subsheaf $\mu^{\ast}L\otimes\mathcal{O}_{\tilde{X}}(-E)$ of $\mu^{\ast}L$ for an effective divisor $E$ on $\tilde{X}$ such that $E+\textrm{except}(\mu)$ has simple normal crossing support. Here $\textrm{except}(\mu)$ is the exceptional divisor of $\mu$. Given $c>0$ we define
\[
\mathscr{I}(f,c|L|)=\mu_{\ast}\mathcal{O}_{\tilde{X}}(K_{\tilde{X}/X}-\lfloor cE\rfloor).
\]
Similarly, $\{\mathscr{I}(f,\frac{c}{k}|L^{k}|)\}_{(k\geqslant0)}$ has a unique maximal element.

\begin{definition}\label{d22}
The relative asymptotic multiplier ideal sheaf associated to $f$, $c$ and $|L|$,
\[
\mathscr{I}(f,c\|L\|)
\]
is defined to be the unique maximal member among the family of ideals $\{\mathscr{I}(f,\frac{c}{k}|L^{k}|)\}$.
\end{definition}

Similarly we can explain the analytic counterpart as we did in the absolute case. Indeed, by definition, $\mathscr{I}(f,c\|L\|)=\mathscr{I}(f,\frac{c}{k}|L^{k}|)$ for some $k$. Let $\rho_{k}$ be the naturally defined homomorphism
\[
\rho_{k}:f^{\ast}f_{\ast}(L^{k})\rightarrow L^{k}.
\]
Let $\mu:\tilde{X}\rightarrow X$ be the smooth modification of $|L^{k}|$ with respect to $f$ such that 
\[
\textrm{Im}(\mu^{\ast}\rho_{k})=\mu^{\ast}L^{k}\otimes\mathcal{O}_{\tilde{X}}(-E).
\] 
Here $E$ is an effective divisor on $\tilde{X}$ such that $E+\textrm{except}(\mu)$ has simple normal crossing support. Consider the ideal sheaf $\mu_{\ast}\mathcal{O}_{\tilde{X}}(-E)$ on $X$. Pick a local coordinate ball $U$ of $Y$, and let $u_{1},...,u_{m}$ be the generators of $\mu_{\ast}\mathcal{O}_{\tilde{X}}(-E)$ on $f^{-1}(U)$. These generators are given by the sections of $\Gamma(f^{-1}(U),L^{k})$ concerning the fact that $\textrm{Im}(\mu^{\ast}\rho_{k})=\mu^{\ast}L^{k}\otimes\mathcal{O}_{\tilde{X}}(-E)$. 

Now let $\varphi_{U}=\log(|u_{1}|^{2}+\cdots+|u_{m}|^{2})$, which is a singular metric on $L^{k}|_{f^{-1}(U)}$. It is then easy to verify that
\[
\mathscr{I}(\frac{c}{k}\varphi_{U})=\mathscr{I}(f,\frac{c}{k}|L^{k}|)\textrm{ on }f^{-1}(U).
\]
Furthermore, if $v_{1},...,v_{m}$ are alternative generators and $\psi_{U}=\log(|v_{1}|^{2}+\cdots+|v_{m}|^{2})$, obviously we have $\mathscr{I}(\frac{c}{k}\varphi_{U})=\mathscr{I}(\frac{c}{k}\psi_{U})$. Hence all the $\mathscr{I}(\frac{c}{k}\varphi_{U})$ patch together to give a globally defined multiplier ideal sheaf $\mathscr{I}(\frac{c}{k}\varphi)$ such that
\[
\mathscr{I}(\frac{c}{k}\varphi)=\mathscr{I}(f,\frac{c}{k}|L^{k}|)=\mathscr{I}(f,c\|L\|)\textrm{ on }X.
\]
One should be careful that $\{f^{-1}(U),\frac{1}{k}\varphi_{U}\}$ won't give a globally defined metric on $L$ in general, since a section of $\Gamma(f^{-1}(U),L^{k})$ is not necessary to extend to the whole space. Hence $\frac{1}{k}\varphi$ is interpreted as the collection of functions $\{f^{-1}(U),\frac{1}{k}\varphi_{U}\}$ by abusing the notation, which is called the collection of (local) singular metrics on $L$ associated to $\mathscr{I}(f,c\|L\|)$. Certainly it is not unique.

\subsection{Abundant divisor}
\label{sec:abundant}
For a nef $\mathbb{Q}$-divisor $D$, its Iitaka dimension $\kappa(D)$ and numerical dimension $\textrm{nd}(D)$ are defined in an obvious way. Then we have the following definition.
\begin{definition}\label{d23}
A nef $\mathbb{Q}$-divisor $D$ is said to be abundant if $\kappa(D)=\textrm{nd}(D)$.
\end{definition}

This notion arises naturally. Moreover, a nef and abundant $\mathbb{Q}$-divisor can be characterised by asymptotic multiplier ideal sheaf as follows.
\begin{theorem}[(Russo, \cite{Rus09})]\label{t21}
Assume that $\kappa(D)\geqslant0$. Then
\[
\mathscr{I}(\|kD\|)=\mathcal{O}_{X}
\]
for all $k$ if and only if $D$ is nef and abundant.
\end{theorem}
\begin{rem}\label{r21}
In his original paper, Russo only talked about the line bundle ($\mathbb{Z}$-divisor, in other words). However, there is no obstacle to extend everything to a $\mathbb{Q}$-divisor with the same argument.
\end{rem}

\subsection{Relative Iitaka dimension}
\label{sec:iitaka}
This part is borrowed from \cite{FEM13}.

Let $f:X\rightarrow Y$ be a surjective morphism between projective manifolds, and $L$ is a line bundle on $X$. Let $l$ be the dimension of a general fibre $F$ of $f$. We have
\begin{proposition}\label{p22}
For every coherent sheaf $\mathcal{G}$ on $X$, there is $C>0$ (independent of $L$) such that
\[
\textrm{rank }f_{\ast}(\mathcal{G}\otimes L^{k})\leqslant Ck^{l} \textrm{ for all }k\gg0.
\]
\begin{proof}
Let us write $L=A\otimes B^{-1}$, with $A$ and $B$ are very ample line bundles. For every $k$, if we choose $E$ general in the complete linear series $|B^{k}|$, then a
local defining function of $E$ is a non-zero divisor on $\mathcal{G}$, in which case we have an inclusion
\[
H^{0}(F,\mathcal{G}\otimes L^{k})\hookrightarrow H^{0}(F,\mathcal{G}\otimes A^{k}).
\]
Since $A$ is very ample, we know that there is a polynomial $P\in\mathbb{Q}[t]$  with $\deg(P)\leqslant l$ such that $h^{0}(F,\mathcal{G}\otimes A^{k})=P(k)$ for $k\gg0$ (see \cite{Laz04}, for example). Therefore $h^{0}(F,\mathcal{G}\otimes L^{k})\leqslant P(k)\leqslant Ck^{l}$ for a suitable $C>0$ and all $k\gg0$.
\end{proof}
\end{proposition}

\begin{definition}\label{d24}
The relative Iitaka dimension $\kappa(L,f)$ of $L$ is the biggest integer $m$ such that there is $C>0$ satisfying
\[
\textrm{rank }f_{\ast}(L^{k})\geqslant Ck^{m} \textrm{ for all }k\gg0
\]
with the convention that $\kappa(L,f)=-\infty$ if $\textrm{rank }f_{\ast}(L^{k})=0$.
\end{definition}
Note that $\kappa(L,f)$ takes value in $\{-\infty,0,1,...,l\}$ by Proposition \ref{p22}. In particular, if $\kappa(L,f)=l$, we say that $L$ is $f$-big.

\section{The vanishing theorem: absolute case}
\label{sec:absolute}
Firstly, let us fix some notations. Let $(X,\Omega)$ be a compact K\"{a}hler manifold of dimension $n$, and let $L$ be a holomorphic line bundle on $X$ endowed with a (singular) Hermitian metric $\varphi\in L^{1}(X)$. More precisely, we first pick a smooth metric $h_{\infty}$ on $L$. Then $\varphi$ defines a new metric $h_{\infty}e^{-\varphi}$ as the weight function. So when we prescribe a singular metric $\varphi$ on $L$, we actually mean that the metric on $L$ is given by $h_{\infty}e^{-\varphi}$. 

Let $A^{p,q}(X,L)$ be the space of smooth $L$-valued $(p,q)$-forms on $X$. Let $\langle\cdot,\cdot\rangle_{\varphi}$ be the pointwise inner product on $A^{p,q}(X,L)$ defined via $h_{\infty}e^{-\varphi}$ and $\Omega$. The $L^{2}$-inner product is then defined by
\[
(\alpha,\beta)_{\varphi}:=\int_{X}\langle\alpha,\beta\rangle_{\varphi}dV_{\Omega}
\]
for $\alpha,\beta\in A^{p,q}(X,L)$, and the norm $\|\cdot\|_{\varphi}$ is induced by $(\cdot,\cdot)_{\varphi}$. In particular, we will directly use the notations $\langle\cdot,\cdot\rangle_{\infty}$, $(\cdot,\cdot)_{\infty}$, $\|\cdot\|_{\infty}$ etc., if the metric involved is $h_{\infty}$ itself (i.e. $\varphi=0$). 
 
\subsection{The Nadel-type vanishing theorem}
\label{sec:nadel}
In this section we will prove Theorem \ref{t13} via an $L^{2}$-estimate developed in \cite{DeP03}.

We firstly do some reductions. Let $\varphi$ be a metric on $L$ that is associated with $\mathscr{I}(\|L\|)$. Then $\varphi$ is of the form that
\[
\varphi=\frac{1}{p}\log(|u_{1}|^{2}+\cdots+|u_{m}|^{2}),
\]
where $\{u_{1},...,u_{m}\}$ is a basis of $H^{0}(X,L^{p})$ for some $p\gg0$ and divisible enough. Let $\mathfrak{a}_{p}$ be the base-ideal of $|L^{p}|$, and let $\mu:\tilde{X}\rightarrow X$ be the smooth modification of $\mathfrak{a}_{p}$ such that
\[
\mu^{\ast}\mathfrak{a}_{p}=\mathcal{O}_{\tilde{X}}(-E),
\]
where $E=\sum\lambda_{i}E_{i}$ has simple normal crossing support. Then locally
\[
\varphi\circ\mu=\frac{1}{p}\sum\lambda_{i}\log|g_{i}|^{2}+\tau,
\]
where $g_{i}$ are local generators of $E_{i}$ and $\tau$ is a (local) smooth quasi-plurisubharmonic function. Hence 
\[
\mathscr{I}(\varphi)=\mu_{\ast}\mathcal{O}_{\tilde{X}}(K_{\tilde{X}/X}-\sum\lfloor\frac{\lambda_{i}}{p}\rfloor E_{i}).
\]
Now, a direct application of Leray's spectral sequence \cite{Har77} implies that
\[
H^{q}(X,K_{X}\otimes L\otimes\mathscr{I}(\varphi))=H^{q}(\tilde{X},K_{\tilde{X}}\otimes\mu^{\ast}L\otimes\mathcal{O}_{\tilde{X}}(-\sum\lfloor\frac{\lambda_{i}}{p}\rfloor E_{i})).
\] 
Let $\hat{L}:=\mu^{\ast}L\otimes\mathcal{O}_{\tilde{X}}(-\sum\lfloor\frac{\lambda_{i}}{p}\rfloor E_{i})$. Then it is enough to consider
\[
H^{q}(\tilde{X},K_{\tilde{X}}\otimes\hat{L}).
\]

Now we can prove Theorem \ref{t13} borrowing the Monge--Amp\`{e}re trick in \cite{DeP03}. We provide the full details for readers' convenience.
\begin{theorem}[(=Theorem \ref{t13})]\label{t31}
Let $X$ be a compact K\"{a}hler manifold of dimension $n$, and let $L$ be a pseudo-effective line bundle on $X$. Then we have
\[
   H^{q}(X,K_{X}\otimes L\otimes\mathscr{I}(\|L\|))=0
\]
for $q>n-\kappa(L)$.
\begin{proof}
Keep the notations. Let $\tilde{L}=\mu^{\ast}L\otimes\mathcal{O}_{\tilde{X}}(-\sum\frac{\lambda_{i}}{p} E_{i})$, which is a $\mathbb{Q}$-bundle. Then 
\[
\hat{L}=\tilde{L}\otimes\mathcal{O}_{\tilde{X}}(\sum\{\frac{\lambda_{i}}{p}\}E_{i}),
\]
where $\{\frac{\lambda_{i}}{p}\}$ refers to the fractional part of $\frac{\lambda_{i}}{p}$. Let $e_{i}$ be the defining section of $E_{i}$. Observe that $\tilde{L}^{p}$, which is a $\mathbb{Z}$-bundle, equipped with the smooth metric $\psi=p\varphi\circ\mu-\sum\lambda_{i}\log|e_{i}|^{2}$, is semi-positive. 

On the other hand, we have
\[
\mu^{\ast}L^{p}=\tilde{L}^{p}\otimes\mathcal{O}_{\tilde{X}}(E).
\]
Then the morphism $\phi_{|\tilde{L}^{p}|}:\tilde{X}\rightarrow\mathbb{P}^{N}$ defined by the free linear series $|\tilde{L}^{p}|$ resolves the indeterminacies of $\phi_{|L^{p}|}$, and in particular it is generically finite over its image. Therefore $\kappa(\tilde{L}^{p})=\kappa(L^{p})=\kappa(L)$. As a result, the numerical dimension $\textrm{nd}(\tilde{L}^{p})$, which is well-defined since $\tilde{L}^{p}$ is semi-positive, is not less than $\kappa(L)$. Hence $c_{1}(\tilde{L}^{p})^{\kappa(L)}\neq0$. 

Fix a K\"{a}hler metric $\omega$ on $\tilde{X}$. Let $\varepsilon>0$. Then $i\Theta_{\tilde{L}^{p},\psi}+\varepsilon\omega$ is a K\"{a}hler form, hence by the Calabi--Yau theorem for complex Monge--Amp\`{e}re equations there exists a smooth metric $\varphi_{\varepsilon}$ on $\tilde{L}^{p}$ such that 
\begin{equation}\label{e31}
i\Theta_{\tilde{L}^{p},\varphi_{\varepsilon}}+\varepsilon\omega>0\quad\textrm{and}\quad(i\Theta_{\tilde{L}^{p},\varphi_{\varepsilon}}+\varepsilon\omega)^{n}=C_{\varepsilon}\omega^{n}.
\end{equation}
Here $C_{\varepsilon}>0$ is a constant such that
\[
C_{\varepsilon}=\frac{\int_{\tilde{X}}(c_{1}(\tilde{L}^{p})+\varepsilon\omega)^{n}}{\int_{\tilde{X}}\omega^{n}}\geqslant C\varepsilon^{n-\kappa(L)}
\]
for a universal constant $C$. Now we consider the metric
\[
\phi_{\varepsilon}=\frac{1}{p}(\delta\varphi_{\varepsilon}+(1-\delta)\psi)+\sum\{\frac{\lambda_{i}}{p}\}\log|e_{i}|^{2}
\]
on $\hat{L}$, where $\delta>0$ is a sufficiently small number which will be fixed later. Since $\varphi_{\varepsilon}$ and $\psi$ are both smooth and $\{\frac{\lambda_{i}}{p}\}<1$, the $L^{2}$-norm $\|\alpha\|^{2}_{\phi_{\varepsilon}}$ is always bounded for any $\alpha\in A^{s,t}(\tilde{X},\hat{L})$. Moreover, by construction,
\[
\begin{split}
i\Theta_{\hat{L},\phi_{\varepsilon}}+\frac{2\varepsilon}{p}\omega&=\frac{1}{p}(\delta(i\Theta_{\tilde{L}^{p},\varphi_{\varepsilon}}+\varepsilon\omega)+(1-\delta)(i\Theta_{\tilde{L}^{p},\psi}+\varepsilon\omega))+\sum\{\frac{\lambda_{i}}{p}\}[E_{i}]+\frac{\varepsilon}{p}\omega \\ 
&\geqslant\frac{\delta}{p}(i\Theta_{\tilde{L}^{p},\varphi_{\varepsilon}}+\varepsilon\omega)+\frac{\varepsilon}{p}\omega.
\end{split}
\]
Here $[E_{i}]$ is the current of integration.

Denote by $0<a_{1}\leqslant\cdots\leqslant a_{n}$ and $0<\hat{a}_{1}\leqslant\cdots\leqslant\hat{a}_{n}$, respectively, the eigenvalues of the curvature forms $i\Theta_{\tilde{L}^{p},\varphi_{\varepsilon}}+\varepsilon\omega$ and $i\Theta_{\hat{L},\phi_{\varepsilon}}+\frac{2\varepsilon}{p}\omega$ at every point $x\in\tilde{X}$, with respect to the base K\"{a}hler metric $\omega(x)$. We find $\hat{a}_{j}\geqslant\frac{\delta}{p}a_{j}+\frac{\varepsilon}{p}$. On the other hand the Monge--Amp\`{e}re equation (\ref{e31}) tells us that
\begin{equation}\label{e32}
a_{1}\cdots a_{n}=C_{\varepsilon}\geqslant C\varepsilon^{n-\kappa(L)}
\end{equation}
everywhere on $\tilde{X}$. 

We apply the Bochner--Kodaira inequality to sections $\alpha$ of type $(n,q)$ with values in $(\hat{L},\phi_{\varepsilon})$. As the curvature eigenvalues of $i\Theta_{\hat{L},\phi_{\varepsilon}}$ are equal to $\hat{a}_{j}-\frac{2\varepsilon}{p}$ by definition, we find
\begin{equation}\label{e33}
\|\bar{\partial}\alpha\|^{2}_{\phi_{\varepsilon}}+\|\bar{\partial}^{\ast}\alpha\|^{2}_{\phi_{\varepsilon}}\geqslant\int_{\tilde{X}}(\hat{a}_{1}+\cdots+\hat{a}_{q}-\frac{2q\varepsilon}{p})|\alpha|^{2}_{\phi_{\varepsilon}}dV_{\omega}.
\end{equation}
Note that (\ref{e33}) is formally true only if $\phi_{\varepsilon}$ is smooth. However, $\phi_{\varepsilon}$ is actually singular along a Zariski closed subset. In this case, it is explained in \cite{DeP03} that in the limit we can recover the same estimates as if we were in the smooth case. We therefore allow ourselves to ignore this minor technical problem from now on.

Now we are ready to prove the desired vanishing result, namely for any $q>n-\kappa(L)$,
\[
H^{q}(X,K_{X}\otimes L\otimes\mathscr{I}(\|L\|))=H^{q}(\tilde{X},K_{\tilde{X}}\otimes\hat{L})=0.
\]
Let us take a cohomology class $[\beta]\in H^{q}(\tilde{X},K_{\tilde{X}}\otimes\hat{L})$. By using the de Rham--Weil isomorphism, we obtain a representative $\beta$ which is a smooth $\hat{L}$-valued $(n,q)$-form. We want to show that $\beta$ is actually a boundary. Since $H^{q}(\tilde{X},K_{\tilde{X}}\otimes\hat{L})$ is a finite dimensional Hausdorff vector space whose topology is induced by the $L^{2}$ Hilbert space topology on the space of forms, it is enough to show that we can approach $\beta$ by $\bar{\partial}$-exact forms in $\|\cdot\|_{\infty}$. Certainly $\|\cdot\|_{\infty}$ here means the $L^{2}$-norm defined by a reference metric $\hat{h}_{\infty}$ on $\hat{L}$ implicitly fixed at the beginning.

As in H\"{o}rmander \cite{Hor65}, we write every form $\alpha$ in the domain of the $L^{2}$-extension (with respect to $\|\cdot\|_{\phi_{\varepsilon}}$) of $\bar{\partial}^{\ast}$ as $\alpha=\alpha_{1}+\alpha_{2}$ with
\[
\alpha_{1}\in\textrm{Ker}\bar{\partial}\quad\textrm{and}\quad\alpha_{2}\in(\textrm{Ker}\bar{\partial})^{\perp}=\overline{\textrm{Im}\bar{\partial}^{\ast}}\subset\textrm{Ker}\bar{\partial}^{\ast}.
\] 
Therefore, since $\beta\in\textrm{Ker}\bar{\partial}$,
\[
\begin{split}
|(\beta,\alpha)_{\phi_{\varepsilon}}|^{2}&=|(\beta,\alpha_{1})_{\phi_{\varepsilon}}|^{2} \\
&\leqslant\int_{\tilde{X}}\frac{1}{\hat{a}_{1}+\cdots+\hat{a}_{q}}|\beta|^{2}_{\phi_{\varepsilon}}dV_{\omega}\int_{\tilde{X}}(\hat{a}_{1}+\cdots+\hat{a}_{q})|\alpha_{1}|^{2}_{\phi_{\varepsilon}}dV_{\omega}.
\end{split}
\]
As $\bar{\partial}\alpha_{1}=0$, an application of (\ref{e33}) to $\alpha_{1}$ shows that the second integral in the right-hand side is bounded above by
\[
\|\bar{\partial}^{\ast}\alpha_{1}\|^{2}_{\phi_{\varepsilon}}+\frac{2q\varepsilon}{p}\|\alpha_{1}\|^{2}_{\phi_{\varepsilon}}\leqslant\|\bar{\partial}^{\ast}\alpha\|^{2}_{\phi_{\varepsilon}}+\frac{2q\varepsilon}{p}\|\alpha\|^{2}_{\phi_{\varepsilon}},
\]
so we finally get
\[
|(\beta,\alpha)_{\phi_{\varepsilon}}|^{2}\leqslant\int_{\tilde{X}}\frac{1}{\hat{a}_{1}+\cdots+\hat{a}_{q}}|\beta|^{2}_{\phi_{\varepsilon}}dV_{\omega}(\|\bar{\partial}^{\ast}\alpha\|^{2}_{\phi_{\varepsilon}}+\frac{2q\varepsilon}{p}\|\alpha\|^{2}_{\phi_{\varepsilon}}).
\]
By the Hahn--Banach theorem, we can find elements $u_{\varepsilon}$, $v_{\varepsilon}$ such that
\[
(\beta,\alpha)_{\phi_{\varepsilon}}=(u_{\varepsilon},\bar{\partial}^{\ast}\alpha)_{\phi_{\varepsilon}}+(v_{\varepsilon},\alpha)_{\phi_{\varepsilon}}
\]
for all $\alpha$. Namely $\beta=\bar{\partial}u_{\varepsilon}+v_{\varepsilon}$, with
\[
\|u_{\varepsilon}\|^{2}_{\phi_{\varepsilon}}+\frac{p}{2q\varepsilon}\|v_{\varepsilon}\|^{2}_{\phi_{\varepsilon}}\leqslant\int_{X}\frac{1}{\hat{a}_{1}+\cdots+\hat{a}_{q}}|\beta|^{2}_{\phi_{\varepsilon}}dV_{\omega}.
\]
As a consequence, the $L^{2}$-distance of $\beta$ to the space of $\bar{\partial}$-exact forms is bounded by $\|v_{\varepsilon}\|_{\phi_{\varepsilon}}$, where
\[
\|v_{\varepsilon}\|^{2}_{\phi_{\varepsilon}}\leqslant\frac{2q\varepsilon}{p}\int_{\tilde{X}}\frac{1}{\hat{a}_{1}+\cdots+\hat{a}_{q}}|\beta|^{2}_{\phi_{\varepsilon}}dV_{\omega}.
\]

We normalise the choice of the potentials $\phi_{\varepsilon}$ so that $\sup_{X}\phi_{\varepsilon}=0$. It is achieved simply by adding suitable constants. From this we infer that
\[
\int_{\tilde{X}}|v_{\varepsilon}|^{2}_{\infty}dV_{\omega}\leqslant\int_{\tilde{X}}|v_{\varepsilon}|^{2}_{\phi_{\varepsilon}}dV_{\omega}\leqslant\frac{2q\varepsilon}{p}\int_{\tilde{X}}\frac{1}{\hat{a}_{1}+\cdots+\hat{a}_{q}}|\beta|^{2}_{\phi_{\varepsilon}}dV_{\omega}.
\]
It remains to show that the right-hand side converges to zero. By construction $\hat{a}_{j}\geqslant\frac{\delta}{p}a_{j}+\frac{\varepsilon}{p}$ and (\ref{e32}) implies
\[
a^{q}_{q}a_{q+1}\cdots a_{n}\geqslant a_{1}\cdots a_{n}\geqslant C\varepsilon^{n-\kappa(L)},
\] 
hence
\[
\frac{1}{a_{1}+\cdots+a_{q}}\leqslant\frac{1}{a_{q}}\leqslant C^{-1/q}\varepsilon^{-(n-\kappa(L))/q}(a_{q+1}\cdots a_{n})^{1/q}.
\]
We infer
\[
\gamma_{\varepsilon}:=\frac{q\varepsilon}{p(\hat{a}_{1}+\cdots+\hat{a}_{q})}\leqslant\min(1,\frac{q\varepsilon}{\delta a_{q}})\leqslant\min(1,C^{\prime}\delta^{-1}\varepsilon^{1-(n-\kappa(L))/q}(a_{q+1}\cdots a_{n})^{1/q}).
\]
We notice that
\[
\int_{\tilde{X}}a_{q+1}\cdots a_{n}dV_{\omega}\leqslant\int_{\tilde{X}}(i\Theta_{\tilde{L}^{p},\varphi_{\varepsilon}}+\varepsilon\omega)^{n-q}\wedge\omega^{q}=(c_{1}(\tilde{L}^{p})+\varepsilon[\omega])^{n-q}[\omega]^{q}\leqslant C^{\prime\prime},
\]
hence the functions $(a_{q+1}\cdots a_{n})^{1/q}$ are uniformly bounded in $L^{1}$-norm as $\varepsilon$ tens to zero. Since $1-\frac{n-\kappa(L)}{q}>0$ by hypothesis, we conclude that $\gamma_{\varepsilon}$ converges almost everywhere to zero as $\varepsilon$ tends to zero.

Now 
\[
\begin{split}
\frac{2q\varepsilon}{p}\int_{\tilde{X}}\frac{1}{\hat{a}_{1}+\cdots+\hat{a}_{q}}|\beta|^{2}_{\phi_{\varepsilon}}dV_{\omega}&=2\int_{\tilde{X}}\gamma_{\varepsilon}|\beta|^{2}_{\infty}e^{-\frac{1}{p}(\delta\varphi_{\varepsilon}+(1-\delta)\psi)}\frac{1}{\Pi|e_{i}|^{2\{\frac{\lambda_{i}}{p}\}}}dV_{\omega} \\
&\leqslant2(\int_{\tilde{X}}e^{-\frac{s\delta}{p}\varphi_{\varepsilon}}dV_{\omega})^{1/s}(\int_{\tilde{X}}\frac{\gamma^{t}_{\varepsilon}|\beta|^{2t}_{\infty}e^{-\frac{t(1-\delta)}{p}\psi}}{\Pi|e_{i}|^{2t\{\frac{\lambda_{i}}{p}\}}}dV_{\omega})^{1/t}
\end{split}
\]
for any $s,t>1$ with $\frac{1}{s}+\frac{1}{t}=1$ by H\"{o}lder's inequality. Then we can take $t=1+\xi$ with $\xi$ small enough such that
\[
\int_{\tilde{X}}\frac{1}{\Pi|e_{i}|^{2t\{\frac{\lambda_{i}}{p}\}}}dV_{\omega}<+\infty.
\]
As $\gamma_{\varepsilon}\leqslant1$, the Lebesgue dominated convergence theorem shows that for every $1>\delta>0$,
\[
\int_{\tilde{X}}\frac{\gamma^{t}_{\varepsilon}|\beta|^{2t}_{\infty}e^{-\frac{t(1-\delta)}{p}\psi}}{\Pi|e_{i}|^{2t\{\frac{\lambda_{i}}{p}\}}}dV_{\omega}
\]
converges to zero as $\varepsilon$ tends to zero. On the other hand, since the curvature forms
\[
i\Theta_{\tilde{L}^{p},\varphi_{\varepsilon}}>-\varepsilon\omega
\]
all sit in $c_{1}(\tilde{L}^{p})$, the family of quasi-plurisubharmonic functions $\{\varphi_{\varepsilon}\}$ is a bounded family (after normalisation) with respect to the $L^{1}$-norm. By standard results of complex potential theory, we conclude that there exits a small constant $\eta>0$ such that $\int_{\tilde{X}}e^{-\eta\varphi_{\varepsilon}}$ is uniformly bounded. (It is carefully explained in \cite{DeP03}.) By choosing $\delta\leqslant\frac{p\eta}{s}$, the integral $\int_{\tilde{X}}e^{-\frac{s\delta}{p}\varphi_{\varepsilon}}$ remains bounded. The proof is then complete.
\end{proof}
\end{theorem}

\subsection{The Kawamata--Veihweg-type vanishing theorem}
\label{sec:kvtype}
In this section we will generalise the Kawamata--Viehweg vanishing theorem using the same Monge--Amp\`{e}re technique.
\begin{theorem}[(=Theorem \ref{t14})]
Let $(X,\omega)$ be a compact K\"{a}hler manifold of dimension $n$, and let $L$ be a holomorphic line bundle on $X$. Assume that $L$ is numerically equivalent to $B+\Delta$, where $B$ is a nef and abundant $\mathbb{Q}$-divisor, and $\Delta=\sum a_{i}\Delta_{i}$ is a $\mathbb{Q}$-divisor with simple normal crossing support and fractional coefficients:
\[
0\leqslant a_{i}<1 \textrm{ for all i}.
\]
Then we have
\[
   H^{q}(X,K_{X}\otimes L)=0
\]
for $q>n-\textrm{nd}(B)$.
\begin{proof}
Firstly we fix some notations. Take a positive integer $p$ such that both $\mathcal{O}_{X}(pB)$ and $\mathcal{O}_{X}(p\Delta)$ are line bundles. Let $\psi$ be the metric on $\mathcal{O}_{X}(pB)$ associated with $\mathscr{I}(P\|B\|)$ where $P>0$ is a sufficiently large number which will be specified later. In particular, $\psi$ depends on $P$ and satisfies that
\[
\mathscr{I}(\frac{P}{p}\psi)=\mathscr{I}(P\|B\|).
\]
On the other hand, $\psi$ is always smooth outside a Zariski closed subset (independent of $P$) and satisfies $i\Theta_{\mathcal{O}_{X}(pB),\psi}\geqslant0$. Since $L$ is numerically equivalent to $B+\Delta$, there exists a flat line bundle $(\tilde{L},\varphi)$ such that $L^{p}=\tilde{L}\otimes\mathcal{O}_{X}(pB)\otimes\mathcal{O}_{X}(p\Delta)$. Let $g_{i}$ be the defining section of $\Delta_{i}$.

Let $\varepsilon>0$. Then $c_{1}(pB)+\varepsilon\omega$ is a K\"{a}hler class, hence by the Calabi--Yau theorem for complex Monge--Amp\`{e}re equations there exists a smooth metric $\varphi_{\varepsilon}$ on $\mathcal{O}_{X}(pB)$ such that 
\begin{equation}\label{e34}
i\Theta_{\mathcal{O}_{X}(pB),\varphi_{\varepsilon}}+\varepsilon\omega>0\quad\textrm{and}\quad(i\Theta_{\mathcal{O}_{X}(pB),\varphi_{\varepsilon}}+\varepsilon\omega)^{n}=C_{\varepsilon}\omega^{n}.
\end{equation}
Here $C_{\varepsilon}>0$ is a constant such that
\[
C_{\varepsilon}=\frac{\int_{X}(c_{1}(\mathcal{O}_{X}(pB))+\varepsilon\omega)^{n}}{\int_{X}\omega^{n}}\geqslant C\varepsilon^{n-\textrm{nd}(B)}
\]
for a universal constant $C$. Now we consider the metric
\[
\phi_{\varepsilon}=\frac{1}{p}(\varphi+\delta\varphi_{\varepsilon}+(1-\delta)\psi)+\sum\log|g_{i}|^{2a_{i}}
\]
on $L$, where $\delta>0$ is a sufficiently small number which will be fixed later. By construction,
\[
\begin{split}
i\Theta_{L,\phi_{\varepsilon}}+\frac{2\varepsilon}{p}\omega&=\frac{1}{p}(\delta(i\Theta_{\mathcal{O}_{X}(pB),\varphi_{\varepsilon}}+\varepsilon\omega)+(1-\delta)(i\Theta_{\mathcal{O}_{X}(pB),\psi}+\varepsilon\omega))+\sum a_{i}[\Delta_{i}]+\frac{\varepsilon}{p}\omega \\ 
&\geqslant\frac{\delta}{p}(i\Theta_{\mathcal{O}_{X}(pB),\varphi_{\varepsilon}}+\varepsilon\omega)+\frac{\varepsilon}{p}\omega.
\end{split}
\]
Here $[\Delta_{i}]$ is the current of integration. 

Denote by $0<a_{1}\leqslant\cdots\leqslant a_{n}$ and $0<\hat{a}_{1}\leqslant\cdots\leqslant\hat{a}_{n}$, respectively, the eigenvalues of the curvature forms $i\Theta_{\mathcal{O}_{X}(pB),\varphi_{\varepsilon}}+\varepsilon\omega$ and $i\Theta_{L,\phi_{\varepsilon}}+\frac{2\varepsilon}{p}\omega$ at every point $x\in X$, with respect to the base K\"{a}hler metric $\omega(x)$. We find $\hat{a}_{j}\geqslant\frac{\delta}{p}a_{j}+\frac{\varepsilon}{p}$. On the other hand the Monge--Amp\`{e}re equation (\ref{e34}) tells us that
\begin{equation}\label{e35}
a_{1}\cdots a_{n}=C_{\varepsilon}\geqslant C\varepsilon^{n-\textrm{nd}(B)}
\end{equation}
everywhere on $X$. 

Let us take a cohomology class $[\beta]\in H^{q}(X,K_{X}\otimes L)$. By using the de Rham--Weil isomorphism, we obtain a representative $\beta$ of the cohomology class which is a smooth $L$-valued $(n,q)$-form. Now it follows the same $L^{2}$-estimate as Theorem \ref{t13} that we can find elements $u_{\varepsilon}$, $v_{\varepsilon}$ such that
\[
\beta=\bar{\partial}u_{\varepsilon}+v_{\varepsilon}
\] 
with
\[
\|u_{\varepsilon}\|^{2}_{\phi_{\varepsilon}}+\frac{p}{2q\varepsilon}\|v_{\varepsilon}\|^{2}_{\phi_{\varepsilon}}\leqslant\int_{X}\frac{1}{\hat{a}_{1}+\cdots+\hat{a}_{q}}|\beta|^{2}_{\phi_{\varepsilon}}dV_{\omega}.
\]
Then through the same estimate we obtain that, when $q>n-\textrm{nd}(B)$, $\gamma_{\varepsilon}:=\frac{q\varepsilon}{p(\hat{a}_{1}+\cdots+\hat{a}_{q})}$ converges almost everywhere to zero as $\varepsilon$ tends to zero.

Now, after a suitable normalisation of the smooth potential $\varphi$, we obtain that
\[
\begin{split}
\frac{2q\varepsilon}{p}\int_{X}\frac{1}{\hat{a}_{1}+\cdots+\hat{a}_{q}}|\beta|^{2}_{\phi_{\varepsilon}}dV_{\omega}&=2\int_{X}\gamma_{\varepsilon}|\beta|^{2}_{\infty}e^{-\frac{1}{p}(\varphi+\delta\varphi_{\varepsilon}+(1-\delta)\psi)}\frac{1}{\Pi|g_{i}|^{2a_{i}}}dV_{\omega} \\
&\leqslant2(\int_{X}e^{-\frac{r\delta}{p}\varphi_{\varepsilon}}dV_{\omega})^{1/r}(\int_{X}e^{-\frac{s(1-\delta)}{p}\psi}dV_{\omega})^{1/s}(\int_{X}\frac{\gamma^{t}_{\varepsilon}|\beta|^{2t}_{\infty}}{\Pi|g_{i}|^{2t a_{i}}}dV_{\omega})^{1/t}
\end{split}
\]
for any $r,s,t>1$ with $\frac{1}{r}+\frac{1}{s}+\frac{1}{t}=1$ by H\"{o}lder's inequality. Then we can take $t=1+\xi$ with $\xi$ small enough such that 
\[
\int_{X}\frac{|\beta|^{2t}_{\infty}}{\Pi|g_{i}|^{2t a_{i}}}dV_{\omega}<+\infty.
\]
As $\gamma_{\varepsilon}\leqslant1$, the Lebesgue dominated convergence theorem shows that 
\[
\int_{\tilde{X}}\frac{\gamma^{t}_{\varepsilon}|\beta|^{2t}_{\infty}}{\Pi|g_{i}|^{2a_{i}}}dV_{\omega}
\]
converges to zero as $\varepsilon$ tends to zero. On the other hand, since the curvature forms
\[
i\Theta_{\mathcal{O}_{X}(pB),\varphi_{\varepsilon}}>-\varepsilon\omega
\]
all sit in $c_{1}(pB)$, the same season as Theorem \ref{t31} implies that there exits a small constant $\eta>0$ such that $\int_{X}e^{-\eta\varphi_{\varepsilon}}$ is uniformly bounded. Fix $r>1$ with $\frac{1}{r}+\frac{1}{t}<1$. By choosing $\delta\leqslant\frac{p\eta}{r}$, the integral $\int_{X}e^{-\frac{r\delta}{p}\varphi_{\varepsilon}}$ remains bounded as $\varepsilon$ tends to zero. At last, we deal with the integral
\[
\int_{X}e^{-\frac{s(1-\delta)}{p}\psi}dV_{\omega}.
\]
As $r$ and $t$ were chosen before, $s$ is computed via the relationship $\frac{1}{r}+\frac{1}{s}+\frac{1}{t}=1$. Since $B$ is nef and abundant, $\mathscr{I}(P\|B\|)=\mathcal{O}_{X}$ for all $P>0$ by Theorem \ref{t21}. Then we take $P$ with $P>s(1-\delta)$. Now $\psi$ is specified accordingly, and satisfies $\int_{X}e^{-\frac{P}{p}\psi}dV_{\omega}<+\infty$. So $\int_{X}e^{-\frac{s(1-\delta)}{p}\psi}dV_{\omega}$ will be bounded, too. In summary, 
\[
\frac{2q\varepsilon}{p}\int_{X}\frac{1}{\hat{a}_{1}+\cdots+\hat{a}_{q}}|\beta|^{2}_{\phi_{\varepsilon}}dV_{\omega}
\]
converges to zero as $\varepsilon$ tends to zero. As a result, $\|v_{\varepsilon}\|_{\phi_{\varepsilon}}$ converges to zero as $\varepsilon$ tends to zero, and $[\beta]=0$ in $H^{q}(X,K_{X}\otimes L)$. The proof is complete.
\end{proof}
\end{theorem}

\section{The vanishing theorem: relative case}
\label{sec:relative}
In this section, we should prove a relative version of the vanishing theorem, i.e. Theorem \ref{t15}. However, if one would like to apply the Monge--Amp\`{e}re technique before, the first thing is to solve a fibrewise Monge--Amp\`{e}re equation, which is not easy. Therefore we prefer to prove it by induction on the dimension, which requires the ambient space to be projective. 

\begin{theorem}[(=Theorem \ref{t15})]\label{t41}
Let $f:X\rightarrow Y$ be a surjective morphism of projective varieties, with $X$ non-singular. Let $L$ be a pseudo-effective line bundle on $X$. Let $l$ be the dimension of a general fibre $F$ of $f$. Then
\[
R^{q}f_{\ast}(K_{X}\otimes L\otimes\mathscr{I}(f,\|L\|))=0
\]
for $q>l-\kappa(L,f)$. 
\begin{proof}
Let $A$ be a very ample line bundle on $X$. Let $H_{1},...,H_{n}$ be generic hypersurfaces in $|A|$, and let
\[
S_{r}:=H_{1}\cap\cdots\cap H_{r}
\]
with the convention that $S=X$ if $r=0$. Then there exists a natural surjective morphism $f_{r}:S_{r}\rightarrow Y_{r}:=f(S_{r})$.

By the standard exact sequence
\[
0\rightarrow K_{X}\otimes L\otimes\mathscr{I}(f,\|L\|)\rightarrow K_{X}\otimes L\otimes A\otimes\mathscr{I}(f,\|L\|)\rightarrow K_{S_{1}}\otimes L|_{S_{1}}\otimes\mathscr{I}(f,\|L\|)|_{S_{1}}\rightarrow0,
\]
we obtain the long exact sequence
\[
\begin{split}
\cdots\rightarrow &R^{q-1}(f_{1})_{\ast}(K_{S_{1}}\otimes L|_{S_{1}}\otimes\mathscr{I}(f,\|L\|)|_{S_{1}})\rightarrow R^{q}f_{\ast}(K_{X}\otimes L\otimes\mathscr{I}(f,\|L\|))\rightarrow \\
&R^{q}f_{\ast}(K_{X}\otimes L\otimes A\otimes\mathscr{I}(f,\|L\|))\rightarrow R^{q}(f_{1})_{\ast}(K_{S_{1}}\otimes L|_{S_{1}}\otimes\mathscr{I}(f,\|L\|)|_{S_{1}})\rightarrow\cdots
\end{split}
\]
Observe that $A$ is ample, $R^{q}f_{\ast}(K_{X}\otimes L\otimes A\otimes\mathscr{I}(f,\|L\|))=0$ for $q>0$ by Nadel's vanishing theorem. (See \cite{Laz04b}, Generalisation 11.2.15, for example.) Hence
\[
R^{q}f_{\ast}(K_{X}\otimes L\otimes\mathscr{I}(f,\|L\|))=R^{q-1}(f_{1})_{\ast}(K_{S_{1}}\otimes L|_{S_{1}}\otimes\mathscr{I}(f,\|L\|)|_{S_{1}}).
\]
On the other hand, let $\varphi$ be the metric associated with $\mathscr{I}(f,\|L\|)$, we have
\[
\mathscr{I}(f,\|L\|)|_{S_{1}}=\mathscr{I}(\varphi)|_{S_{1}}=\mathscr{I}(\varphi|_{S_{1}}).
\]
The last equality is due the Ohsawa--Takegoshi extension theorem and Fubini's theorem. One could refer to \cite{FuM21} for a complete discussion. We then obtain that
\[
R^{q}f_{\ast}(K_{X}\otimes L\otimes\mathscr{I}(f,\|L\|))=R^{q-1}(f_{1})_{\ast}(K_{S_{1}}\otimes L|_{S_{1}}\otimes\mathscr{I}(\varphi|_{S_{1}})).
\]
Repeating this procedure, we will eventually obtain that
\[
R^{q}f_{\ast}(K_{X}\otimes L\otimes\mathscr{I}(f,\|L\|))=R^{1}(f_{q-1})_{\ast}(K_{S_{q-1}}\otimes L|_{S_{q-1}}\otimes\mathscr{I}(\varphi|_{S_{q-1}})).
\]

Since the relative Iitaka dimension is not decreasing in restriction, $L|_{S_{q-1}}$ is actually $(f_{q-1})$-big when $q>l-\kappa(L,f)$. In this situation we will prove that 
\[
R^{1}(f_{q-1})_{\ast}(K_{S_{q-1}}\otimes L|_{S_{q-1}}\otimes\mathscr{I}(\varphi|_{S_{q-1}}))=0,
\]
which leads to the desired vanishing result.

Fix $p\gg0$ and divisible enough that computes $\mathscr{I}(f,\|L\|)$. Let $\mathfrak{a}_{p,f}$ be the base-ideal of $|L^{p}|$ relative to $f$. By definition, it is the image of the natural morphism
\[
f^{\ast}f_{\ast}(L^{p})\otimes L^{-p}\rightarrow\mathcal{O}_{X}.
\]
Let $\mu:\tilde{S}_{q-1}\rightarrow S_{q-1}$ be the smooth modification of $\mathfrak{a}_{p,f}|_{S_{q-1}}$ such that
\[
\mu^{\ast}(\mathfrak{a}_{p,f}|_{S_{q-1}})=\mathcal{O}_{\tilde{S}_{q-1}}(-E),
\]
where $E=\sum\lambda_{i}E_{i}$ has simple normal crossing support. Then locally
\[
\varphi|_{S_{q-1}}\circ\mu=\frac{1}{p}\sum\lambda_{i}\log|g_{i}|^{2}+\tau,
\]
where $g_{i}$ are local generators of $E_{i}$ and $\tau$ is a (local) smooth quasi-plurisubharmonic function. Hence 
\[
\mathscr{I}(\varphi|_{S_{q-1}})=\mu_{\ast}\mathcal{O}_{\tilde{S}_{q-1}}(K_{\tilde{S}_{q-1}/S_{q-1}}-\sum\lfloor\frac{\lambda_{i}}{p}\rfloor E_{i}).
\]
Now, a direct application of Leray's spectral sequence implies that
\[
R^{1}(f_{q-1})_{\ast}(K_{S_{q-1}}\otimes L|_{S_{q-1}}\otimes\mathscr{I}(\varphi|_{S_{q-1}}))=R^{1}(f_{q-1}\circ\mu)_{\ast}(K_{\tilde{S}_{q-1}}\otimes\mu^{\ast}L|_{S_{q-1}}\otimes\mathcal{O}_{\tilde{S}_{q-1}}(-\sum\lfloor\frac{\lambda_{i}}{p}\rfloor E_{i})).
\] 
Let $\tilde{L}=\mu^{\ast}L|_{S_{q-1}}\otimes\mathcal{O}_{\tilde{S}_{q-1}}(-\sum\frac{\lambda_{i}}{p} E_{i})$, which is a $\mathbb{Q}$-bundle. Then 
\[
\mu^{\ast}L|_{S_{q-1}}\otimes\mathcal{O}_{\tilde{S}_{q-1}}(-\sum\lfloor\frac{\lambda_{i}}{p}\rfloor E_{i})=\tilde{L}\otimes\mathcal{O}_{\tilde{S}_{q-1}}(\sum\{\frac{\lambda_{i}}{p}\}E_{i}),
\]
where $\{\frac{\lambda_{i}}{p}\}$ refers to the fractional part of $\frac{\lambda_{i}}{p}$. 

Let $\{U_{\alpha}\}$ be an affine covering of $Y$. We write $\varphi$ as $\varphi=\{f^{-1}(U_{\alpha}),\varphi_{\alpha}\}$ by abusing the notation (see Sect. \ref{sec:asymptotic}). Let $e_{i}$ be the defining section of $E_{i}$. Observe that $\tilde{L}^{p}$, which is a $\mathbb{Z}$-bundle, equipped with 
\[
\psi=p\varphi_{\alpha}|_{S_{q-1}}\circ\mu-\sum\lambda_{i}\log|e_{i}|^{2}
\] 
is semi-positive. Hence $\tilde{L}$ is nef on every $\mu^{-1}(f^{-1}_{q-1}(U_{\alpha}))$. Equivalently, $\tilde{L}$ is $(f_{q-1}\circ\mu)$-nef. On the other hand, we have
\[
\mu^{\ast}L^{p}|_{S_{q-1}}=\tilde{L}^{p}\otimes\mathcal{O}_{\tilde{S}_{q-1}}(E).
\]
Then the morphism $\phi_{|\tilde{L}^{p}|/Y}:\tilde{S}_{q-1}\rightarrow\mathbb{P}^{N}$ over $Y$ defined by the $(f_{q-1}\circ\mu)$-free relative linear series \cite{Laz04} of $\tilde{L}^{p}$ resolves the indeterminacies of $\phi_{|L^{p}|_{S_{q-1}}|/Y}$, and in particular it is generically finite over its image. Therefore $\kappa(\tilde{L}^{p},f_{q-1}\circ\mu)=\kappa(L^{p}|_{S_{q-1}},f_{q-1})=\kappa(L|_{S_{q-1}},f_{q-1})$. In summary, $\tilde{L}$ is $(f_{q-1}\circ\mu)$-nef and $(f_{q-1}\circ\mu)$-big.

Now applying the relative version of the Kawamata--Viehweg vanishing theorem (See \cite{Laz04b}, Generalisation 9.1.22, for example), we eventually obtain that
\[
0=R^{1}(f_{q-1}\circ\mu)_{\ast}(K_{\tilde{S}_{q-1}}\otimes\lceil\tilde{L}\rceil)=R^{1}(f_{q-1}\circ\mu)_{\ast}(K_{\tilde{S}_{q-1}}\otimes\mu^{\ast}L|_{S_{q-1}}\otimes\mathcal{O}_{\tilde{S}_{q-1}}(-\sum\lfloor\frac{\lambda_{i}}{p}\rfloor E_{i})).
\]
Here $\lceil\tilde{L}\rceil$ means the round-up. The proof is complete.
\end{proof}
\end{theorem}

\end{document}